\documentclass[10pt
]{amsart}
\usepackage{amssymb,amscd, epsfig, mathrsfs, xypic, amsmath,amsthm, color, dsfont}
\usepackage{hyperref}
\usepackage{blindtext}
\usepackage{enumerate}

\xyoption{all}
\newtheorem{thm}{Theorem}[section]

\newtheorem{df-pr}[thm]{Definition-Proposition}

\theoremstyle{definition}

\newcommand{\ZZ}{{\mathbb Z}}

\newsavebox{\savepar}

\numberwithin{equation}{section}

\newcounter{labelflag} 
\newcommand{\labelon}{\setcounter{labelflag}{1}}
\newcommand{\Label}[1]{\ifnum\thelabelflag=1\ifmmode
\makebox[0in][l]{\qquad\fbox{\rm#1}} \else
\marginpar{\vspace{0.7\baselineskip} \hspace{-1.1\textwidth}
\fbox{\rm#1}} \fi \fi \label{#1} } \labelon

\begin{document} 
\title{An algebraic proof of determinant formulas of Grothendieck polynomials}
\author{Tomoo Matsumura}
\maketitle 

\begin{abstract}
We give an algebraic proof of the determinant formulas for factorial Grothendieck polynomials obtained by Hudson--Ikeda--Matsumura--Naruse in \cite{HIMN} and by Hudson--Matsumura in \cite{HudsonMatsumura}.
\end{abstract}
\section{Definition and Theorem}
In \cite{LascouxSchutzenberger3} and \cite{Lascoux1}, Lascoux and Sch{\"u}tzenberger introduced \emph{(double) Grothendieck polynomials} indexed by permutations as representatives of $K$-theory classes of structure sheaves of Schubert varieties in a full flag variety. In \cite{GrothendieckFomin} and \cite{DoubleGrothendieckFomin}, Fomin and Kirillov described them combinatorially in the framework of Yang-Baxter equations. Let $x=(x_1,\dots,x_d)$, $b=(b_1,b_2,\dots)$ be sets of indetermiants. A \emph{Grassmannian permutation} with descent at $d$ corresponds to a \emph{partition} $\lambda$ of length at most $d$, {\it i.e.} a sequence of non-negative integers $\lambda=(\lambda_1,\dots,\lambda_d)$ such that $\lambda_i\geq \lambda_{i+1}$ for each $i=1,\dots,d-1$. For such permutation, Buch \cite{BuchLRrule} gave a combinatorial expression of the corresponding Grothendieck polynomial $G_{\lambda}(x)$ as a generating series of \emph{set-valued tableaux}, a generalization of semi-standard Young tableaux by allowing a filling of a box in the Young diagram to be a set of integers. In \cite{McNamara}, McNamara gave an expression of \emph{factorial} (double) Grothendieck polynomials $G_{\lambda}(x|b)$ also in terms of set-valued tableaux. 

In this paper, we prove the following Jacobi--Trudi type determinant formula for $G_{\lambda}(x|b)$. For each non-negative integer $k$ and an integer $m$, let $G_{m}^{(k)}(x|b)$ be a function of $x$ and $b$ given by 
\[
G^{(k)}(u):=\sum_{m\in \ZZ} G_m^{(k)}(x|b) u^m := \frac{1}{1+\beta u^{-1}}\prod_{i=1}^d\frac{1+\beta x_i}{1- x_i u} \prod_{j=1}^{k} (1+  (u+\beta)b_j),
\]
where $\beta$ is a formal variable of degree $-1$ and $\displaystyle\frac{1}{1+\beta u^{-1}}$ is expanded as $\sum_{s\geq 0} (-1)^s\beta^s u^{-s}$. We use the generalized binomial coefficients $\binom{n}{i}$ given by $(1+x)^{n} = \sum_{i\geq 0} \binom{n}{i} x^i$ for $n\in \ZZ$ with the convention that $\binom{n}{-i}=0$ for all integers $i>0$.
\begin{thm}\label{mainthm}
For each partition $\lambda$ of length at most $d$, we have
\begin{eqnarray}
G_{\lambda}(x|b) 
&=& \det\left(\sum_{s\geq 0} \binom{i-d}{s} \beta^sG_{\lambda_i+j-i+s}^{(\lambda_i+d-i)}(x|b)\right)_{1\leq i,j\leq d}\label{HM2}\\
&=& \det\left(\sum_{s\geq 0} \binom{i-j}{s} \beta^sG_{\lambda_i+j-i+s}^{(\lambda_i+d-i)}(x|b)\right)_{1\leq i,j\leq d}\label{HIMN1}.
\end{eqnarray}
In particular, we have $G_{(k,0,\dots,0)}(x|b) = G_k^{(k+d-1)}(x|b)$.
\end{thm}
The formulas (\ref{HM2}) and (\ref{HIMN1}) were originally obtained in the context of degeneracy loci formulas for flag bundles by Hudson--Matsumura in \cite{HudsonMatsumura} and Hudson--Ikeda--Matsumura--Naruse in \cite{HIMN} respectively. The proof in this paper is purely algebraic, generalizing Macdonald's argument in \cite[(3.6)]{Macdonald} for Jacobi--Trudi formula of Schur polynomials. We use the following ``bi-alternant'' formula of $G_{\lambda}(x|b)$ described by Ikeda--Naruse in \cite{IkedaNaruse}:
\begin{equation}\label{eqbialt}
G_{\lambda}(x|b) =\frac{\det\left([x_j|b]^{\lambda_i+d-i}(1+\beta x_j)^{i-1}\right)_{1\leq i,j\leq d}}{\prod_{1\leq i<j\leq d}(x_i-x_j)}.
\end{equation}
Here we denote $x\oplus y := x+y+\beta xy$ and $[y|b]^k := (y \oplus b_1)\cdots (y \oplus b_k)$ for any variable $x,y$. Note that the Grothendieck polynomial $G_{\lambda}(x)$ described in \cite{BuchLRrule} coincides with $G_{\lambda}(x|b)$ by setting $\beta=-1$ and $b_i=0$.

In the non-factorial case, Jacobi--Trudi type formulas different from the ones in Theorem \ref{mainthm} have been also obtained by Lenart in \cite{Lenart2000} and by Kirillov in \cite{Kirillov2015} (see also the papers \cite{LascouxNaruse} by Lascoux--Naruse and \cite{Yeliussizov} by Yeliussizov). It is worth mentioning that in \cite{BuchLRrule} Buch  obtained the Littlewood-Richardson rule for the structure constants of Grothendieck polynomials $G_{\lambda}(x)$, and hence the Schubert structure constants of the $K$-theory of Grassmannians (see also the paper \cite{IkedaShimazaki} by Ikeda-Shimazaki  for another proof). For the equivariant $K$-theory of Grassmannians (or equivalently for the factorial Grothendieck polynomials $G_{\lambda}(x|b)$), the structure constants were recently determined by Pechenik and Yong in \cite{PechenikYongKTG} by introducing a new combinatorial object called \emph{genomic tableaux}. 

\section{Proof of (\ref{HM2})}
By (\ref{eqbialt}), it suffices to show the identity
\[
\frac{\det\left([x_j|b]^{a_i+d-i}(1+\beta x_j)^{i-1}\right)_{1\leq i,j\leq d}}{\prod_{1\leq i<j\leq d}(x_i-x_j)} 
= \det\left(\sum_{s\geq 0} \binom{i-d}{s} \beta^sG_{a_i+j-i+s}^{(a_i+d-i)}(x|b)\right)_{1\leq i,j\leq d},
\]
for each $(a_1,\dots,a_d)\in \ZZ^d$ such that $a_i+d-i\geq 0$.
For each integer $j$ such that $1\leq j \leq d$, we let
\[
E^{(j)}(u):=\sum_{p=0}^{d-1} e_p^{(j)}(x) u^p := \prod_{1\leq i\leq d \atop{i\not=j}} (1+x_i u).
\]
We denote $\bar y:=\displaystyle\frac{-y}{1+\beta y}$. Since $1+ (u+\beta)y=\displaystyle\frac{1- \bar yu}{1+\beta \bar y}$, we have
\[
G^{(k)}(u)= \frac{1}{1+\beta u^{-1}}\prod_{i=1}^d\frac{1+\beta x_i}{1- x_i u} \prod_{\ell=1}^{k} \frac{1-\bar b_\ell u}{1+\beta \bar b_\ell}.
\]
Consider the identity
\begin{equation}\label{eq1}
G^{(k)}(u)E^{(j)}(-u)
= \frac{1}{1+\beta u^{-1}}\frac{1}{1- x_j u}\prod_{i=1}^d(1+\beta x_i) \prod_{\ell=1}^{k} \frac{1-{\bar b}_\ell u}{1+\beta {\bar b}_\ell}.
\end{equation}
\begin{eqnarray*}
\sum_{p=0}^{d-1}G_{m-p}^{(k)}(x|b)  (-1)^pe_p^{(j)}(x)
&=& x_j^{m-k}  \frac{\prod_{\ell=1}^k  (x_j-{\bar b}_{\ell})}{\prod_{\ell=1}^{k} (1+\beta {\bar b}_\ell)} \prod_{1\leq i\leq d\atop{i\not=j}}(1+\beta x_i).
\end{eqnarray*}
Since $\displaystyle\frac{y- \bar b}{1+\beta \bar b} = y \oplus b$, we have 
\begin{eqnarray}\label{eq2}
\sum_{p=0}^{d-1}G_{m-p}^{(k)}(x|b)  (-1)^pe_p^{(j)}(x)
&=& x_j^{m-k}[x_j|b]^k\prod_{1\leq i\leq d \atop{i\not=j}}(1+\beta x_i), \ \ \ \ (m\geq k).
\end{eqnarray}

Consider the matrices
\begin{eqnarray*}
H:=\left(\sum_{s\geq 0} \binom{i-{d}}{s}\beta^s G_{a_i+j-i+s}^{(a_i+d-i)}(x|b)\right)_{1\leq i,j\leq d} \ \ \ \mbox{and} \ \ \ M:=\left((-1)^{d-i}  e_{d-i}^{(j)}(x)\right)_{1\leq i,j\leq d}.
\end{eqnarray*}
By using (\ref{eq2}), we find that the $(i,j)$-entry of $HM$ is
\begin{eqnarray*}
(HM)_{ij}
&=&[x_j|b]^{a_i+d-i}\left(1+\beta x_j\right)^{i-d-1}\prod_{1\leq t\leq d}(1+\beta x_t).
\end{eqnarray*}
By taking the determinant of $HM$, the factor $\prod_{1\leq j\leq d}(1+\beta x_j)^{-d}\prod_{1\leq t\leq d}(1+\beta x_t)^d$ which turns to be $1$ comes out, and therefore we obtain
\[
\det H \det M = \det \left([x_j| b]^{a_i+d-i}\left(1+\beta x_j\right)^{i-1} \right)_{1\leq i, j\leq d}.
\] 
By dividing by $\det M$, we obtain the desired identity since $\det M = \prod_{1\leq i<j\leq d} (x_i-x_j)$ (see \cite[p.42]{Macdonald}). \qed
\section{Proof of (\ref{HIMN1})}
By (\ref{eqbialt}), it suffices to show the identity
\[
\frac{\det\left([x_j|b]^{a_i+d-i}(1+\beta x_j)^{i-1}\right)_{1\leq i,j\leq d}}{\prod_{1\leq i<j\leq d}(x_i-x_j)} 
= \det\left(\sum_{s\geq 0} \binom{i-j}{s} \beta^sG_{a_i+j-i+s}^{(a_i+d-i)}(x|b)\right)_{1\leq i,j\leq d}.
\]
for each $(a_1,\dots,a_d)\in \ZZ^d$ such that $a_i+d-i\geq 0$. Fix $j$ such that $1\leq j \leq d$. Let
\[
{\overline{E}}^{(j)}(u):=\sum_{p=0}^{d-1} e_p^{(j)}(-\bar x) u^p := \prod_{1\leq i\leq d \atop{i\not=j}} (1-\bar x_i u).
\]
Since $1+ (u+\beta)y=\displaystyle\frac{1- \bar yu}{1+\beta \bar y}$, we have the identity
\begin{eqnarray}\label{eq100}
G^{(k)}(u){\overline{E}}^{(j)}(-u-\beta)
&=&\frac{1}{1+\beta u^{-1}}\frac{1+\beta x_j}{1- x_j u}\prod_{1\leq \ell\leq k} \frac{1-\bar b_{\ell} u}{1+\beta \bar b_{\ell}}.
\end{eqnarray}
By comparing the coefficient of $u^m, m\geq k$ in (\ref{eq100}) we obtain
\begin{eqnarray}\label{eq200}
\sum_{p=0}^{d-1}\sum_{s=0}^p \binom{p}{s}\beta^s G_{m-p+s}^{(k)}(x|b)  (-1)^pe_p^{(j)}(-\bar x)
&=&  x_j^{m-k} \prod_{1\leq \ell\leq k} \frac{x_j-\bar b_{\ell}}{1+\beta \bar b_{\ell}} = x_j^{m-k} [x_j| b]^k
\end{eqnarray}
where the last equality follows from the identity $\displaystyle\frac{x-\bar y}{1+\beta \bar y} = x\oplus y$ for any variable $x,y$.
Consider the matrices
\begin{eqnarray*}
H':=\left(\sum_{s\geq 0} \binom{i-j}{s}\beta^s G_{a_i+j-i+s}^{(a_i+d-i)}(x|b)\right)_{1\leq i,j\leq d} \ \ \ \mbox{and} \ \ \ \overline{M}:=\left((-1)^{d-i}  e_{d-i}^{(j)}(-\bar x)\right)_{1\leq i,j\leq d}.
\end{eqnarray*}
We write the $(i,j)$-entry of the product $H'\overline{M}$ as
\begin{eqnarray*}
(H'\overline{M})_{ij}
&=&\sum_{p=0}^{d-1}\sum_{s\geq 0} \binom{i-{d+p}}{s}\beta^s G_{a_i+d-i+s-p}^{(a_i+d-i)}(x|b)(-1)^{p}  e_{p}^{(j)}(-\bar x).
\end{eqnarray*}
By writing $\displaystyle\binom{i-{d+p}}{s} = \sum_{\ell\geq 0}  \binom{i-d}{\ell}\binom{p}{s-\ell}$ using a well-known identity of binomial coefficients and then applying (\ref{eq200}), we obtain
\begin{eqnarray*}
(H'\overline{M})_{ij}
&=&[x_j| b]^{a_i+d-i}(1+\beta x_j)^{i-1}(1+\beta x_j)^{1-d}.
\end{eqnarray*}
By taking the determinant of $H'\overline{M}$, we have
\[
\det H'\det\overline{M} =\left(\prod_{1\leq j\leq d} (1+\beta x_j)^{1-d} \right)\det\left([x_j| b]^{a_i+d-i}(1+\beta x_j)^{i-1}\right)_{1\leq i,j\leq d}.
\]
Since we have (see \cite[p.42]{Macdonald})
\[
\det \overline{M} = \prod_{1\leq i<j\leq d} (\bar x_j - \bar x_i) = \prod_{1\leq i<j\leq d} \frac{x_i-x_j}{(1+\beta x_i)(1+\beta x_j)} = \prod_{1\leq i\leq d} (1+\beta x_i)^{1-d} \prod_{1\leq i<j\leq d} (x_i-x_j),
\]
we obtain the desired identity. \qed

\

\

\noindent\textbf{Acknowledgements.} 
The author would like to thank Takeshi Ikeda for useful discussions. The author is supported by 	Grant-in-Aid for Young Scientists (B) 16K17584.
\bibliography{references}{}
\bibliographystyle{acm}
\end{document}